\documentstyle[12pt,righttag,intlim]{amsart}
\title{Heat invariants of Riemannian manifolds}
\author{Iosif Polterovich}
\address{Department of Theoretical Mathematics, The Weizmann Institute
of Science, Rehovot 76100, Israel}
\email{iossif@@wisdom.weizmann.ac.il}
\def \phi{\varphi}
\def \epsilon{\varepsilon}
\numberwithin{equation}{subsection}
\theoremstyle{definition}

\theoremstyle{plain}
\newtheorem{lemma}[equation]{Lemma}
\newtheorem{theorem}[equation]{Theorem}

\begin{document}
\maketitle
\begin{abstract}
We calculate heat invariants of arbitrary Riemannian manifolds without 
boundary. Every heat invariant is expressed in terms of powers of the 
Laplacian and the distance function. Our approach is based on a 
multi-dimensional generalization of the Agmon-Kannai method.
An application to computation of the Korteweg-de Vries
hierarchy is also presented.
\end{abstract}
\section{Introduction and main results}
\subsection{Heat invariants}
Let $M$ be a $d$-dimensional Riemannian manifold 
without boundary with a  metric $(g_{ij})$, and
$\Delta$ be the {\it Laplace-Beltrami} 
operator (or simply the {\it Laplacian}) on $M$.
In local coordinates $(x_1,\dots,x_d)$ the Laplacian is given by
\begin{equation}
\Delta f = - \frac{1}{\sqrt{g}}\sum_{i,j=1}^d 
\frac{\partial(\sqrt{g} g^{ij}(\partial f/\partial x_i))}{\partial x^j},
\end{equation}
where $g=$det$(g_{ij})$, and $(g^{ij})$ denotes the inverse of the 
matrix $(g_{ij})$.

The {\it heat kernel} $K(t,x,y)$ is the fundamental solution of the 
heat equation 
$$(\frac{\partial}{\partial t}+\Delta)f=0.$$ 
The function $K(t,x,y)$
is analytic in $t>0$ and $C^{\infty}$ in 
$x$ and $y$, and has 
the following asymptotic expansion on the diagonal as $t\to 0+$ (see [G1]):
$${\mbox K}(t,x,x) \sim \sum_{n=0}^{\infty} a_n(x)t^{n-\frac{d}{2}}
$$
It is called the Minakshisundaram-Plejel asymptotic expansion ([MP]).
The coefficients $a_n(x)$ are (local) {\it heat invariants} 
of the manifold $M$. They are homogeneous polynomials of degree $2n$ in the 
derivatives of the Riemannian metric $\{g^{ij}\}$ at the point $x$ (see [G2]).
Integrating $a_n(x)$ over
the manifold one gets the coefficients $a_n$ of the expansion for 
the trace  of the heat operator $e^{-t\Delta}$:  
\begin{equation}
\label{tr}
\sum_i e^{-t\lambda_i} \sim  \sum_{n=0}^{\infty}\left(\,\,\int\limits_M 
a_n(x)\sqrt{g}dx \right) t^{n-\frac{d}{2}} \sim  \sum_{n=0}^{\infty}a_n
t^{n-\frac{d}{2}}.
\end{equation}

Computation of heat invariants is a well-known problem in spectral geometry
(see [Be], [G1], [Ch], [G3], [Ro]) which has various 
applications ([F], [P2]). The first method for derivation of heat kernel
asymptotics is due to Seeley ([Se]). This method was deveoped later 
by Gilkey (see Theorem 1.3 in [G1]) who presented a way to get recursive 
formulas for the heat invariants. 

However, explicit formulas for $a_n(x)$ in arbitrary dimension 
existed only for $n\le 5$ ([MS], [Sa], [Av], [vdV]). 
The reason for this is the combinatorial complexity
of $a_n(x)$ which is increasing very rapidly with the growth of $n$. 
For the higher heat invariants only partial information is known 
([BG\O], [OPS]). Let us  also mention interesting recursive formulas for 
$a_n(x)$ obtained in [Xu].

In this  paper we calculate all heat invariants $a_n(x)$ for an arbitrary 
Riemannian manifold without boundary in terms of
powers of the Laplacian and the distance function. 
Note that the heat invariants were initially given by a recursive system 
of differential equations involving exactly the Laplacian and the distance 
function  (see [MP]).

\subsection{Main result}
Given a point $x\in M$ denote by $\rho_x:M\to {\Bbb R}$ 
the corresponding distance function: 
for every $y\in M$ the distance between 
the points $y$ and $x$ is  $\rho_x(y)$.

\begin{theorem}
\label{most}
Heat invariants $a_n(x)$ are equal to 
\begin{equation*}
a_n(x)=(4\pi)^{-d/2}(-1)^n\sum_{j=0}^{3n}
\binom{3n+\frac{d}{2}}{j+\frac{d}{2}}
\frac{1}{4^j \, j! \, (j+n)!}\left.
\Delta^{j+n}(\rho_x(y)^{2j})\right|_{y=x}.
\end{equation*}
\end{theorem}

The binomial coefficients for $d$ odd are defined by (\ref{bc}).

\subsection{Structure of the paper}
In [P1], [P2] we have developed a method  for computation of heat
invariants based on the Agmon-Kannai asymptotic expansion 
of resolvent kernels of elliptic operators ([AK]). 
In [P1] this method is used to obtain explicit
formulas for the heat invariants of $2$-dimensional Riemannian 
mainfolds, in [P2] --- for computation of the Korteweg-de Vries (KdV) 
hierarchy via heat kernel coefficients of the 
$1$-dimensional Schr\"odinger operator. 
In this paper we present a 
multi-dimensional generalization of the Agmon-Kannai method  which 
is described in section 2.3. 
In section 3.1 we apply it to get 
formulas for the heat invariants in normal coordinates. It turns out that 
combinatorial coefficients in these formulas can be substantially simplified
which is done in section 4.1. In section 4.2 we present $a_n(x)$ in a 
completely invariant form and prove Theorem \ref{most}. 
Main result allows to simplify the formulas for 
the KdV hierarchy obtained in [P2]. This is shown in sections 5.1 and 5.2.

\section{Asymptotics of derivatives of the resolvent}
\subsection{A modification of Agmon-Kannai expansion}
The original Agmon--Kan\-nai theorem ([AK]) deals with asymptotic behaviour
of resolvent kernels of elliptic operators. In [P1] we have obtained
a concise reformulation of this theorem which is suitable
for computation of heat invariants. 
We start with some notations.

Let $H$ be a a self--adjoint elliptic differential 
operator of order $p$ on a Riemannian manifold $(M, g_{ij})$ 
of dimension $d<p$ and let $H_0$ be the operator obtained by
freezing the coefficients of the principal part $H'$ of the operator
$H$ at some point $x\in M$: $H_0=H'(x)$.  
Denote by $R_\lambda(x,y)$ the kernel of the resolvent
$R_\lambda=(H-\lambda)^{-1}$, and by $F_{\lambda}(x,y)$ --- the 
kernel of $F_{\lambda}=(H_0-\lambda)^{-1}$.
\begin{theorem}$\operatorname{([P1])}.$
\label{ak}
The resolvent kernel $R_{\lambda}(x,y)$ 
has the following asymptotic representation on the diagonal
as $\lambda \to \infty$:
\begin{equation}
\label{s}
R_\lambda(x,x) \sim \frac{1}{\sqrt{g}}\sum_{m=0}^\infty 
X_m F_\lambda^{m+1}(x,x),
\end{equation}
where the operators $X_m$ are defined by:
\begin{equation}
\label{y}
X_m=\sum_{k=0}^m (-1)^k \binom{m}{k} H^k H_0^{m-k}, \,\, m\ge 0.
\end{equation}
\end{theorem}
\subsection{Derivatives of the resolvent}
The main obstruction for using Theorem \ref{ak} directly for
computation of heat invariants of a $d$-dimensional Riemannian manifold
is the condition $d<p$, where $p=2$ is the order of the Laplacian.
In [P1] we avoid this difficulty for 2-dimensional manifolds
taking the difference of resolvents. However, in  the general case one
should consider derivatives of the resolvent kernel (cf. [AvB]).
\begin{lemma}
The following asymptotic expansion on the diagonal holds for the derivatives 
of the resolvent kernel of the Laplacian on a $d$-dimensional Riemannian
manifold $M$:
\begin{equation}
\label{res}
\frac{d^s}{d\lambda^s}R_\lambda(x,x)\sim \sum_{n=0}^{\infty}
\Gamma(s+n-\frac{d}{2}+1)a_n(x) (-\lambda)^{\frac{d}{2}-s-n-1},
\,\,\, s\ge d/2,
\end{equation}
where $a_n(x)$ are heat invariants of the manifold $M$.
\end{lemma}
\noindent {\bf Proof.}
Let $\operatorname{Re} \lambda <0$. We have (formally):
$$\int_0^\infty e^{-t(\Delta-\lambda)}dt=\frac{1}{\Delta-\lambda}$$
Differentiating $R_{\lambda}$ 
$s$ times with respect to $\lambda$ we get
a self--adjoint operator from 
$L^2(M)$ into the Sobolev space ${\mbox H}^{2s+2}(M)$.
Since $2s+2> \operatorname{dim} M$ this operator has a continuous kernel 
(see [AK]).
Taking into account (\ref{tr}) we formally have:
\begin{equation}
\label{ros}
\frac{d^s}{d\lambda^s}\left(\frac{1}{\Delta-\lambda}\right)=
\int_0^\infty t^se^{- t(\Delta-\lambda)}dt
\sim \sum_{n=0}^\infty a_n
\int_0^\infty 
t^{s+n-d/2}e^{\lambda t}dt
\end{equation}

The asymptotic expansion in (\ref{ros}) is obviously valid
if we integrate  over a finite interval $[0,T]$. 
In order to  show that it remains true in our case as well we need an
additional argument.

Indeed, it is well-known (for example, see [Da]) that
$$|e^{-t\Delta}|\le ct^{-d/2}.$$

Therefore we have:
$$
\left|\int_0^\infty t^se^{- t(\Delta-\lambda)}dt-
\int_0^T t^se^{- t(\Delta-\lambda)}dt\right|\le
\int_T^\infty ct^{s-d/2}e^{\lambda t}dt.$$
Let us estimate the second integral. We have:
$$
\int_T^\infty t^{s-d/2}e^{\lambda t}dt\le e^{-\epsilon T}\int_T^{\infty}
t^{s-d/2}e^{(\lambda+\epsilon)t}dt\le e^{-\epsilon T}\frac{\Gamma(s-d/2+1)}
{(\lambda+\epsilon)^{s-d/2+1}}$$
Take $\epsilon=\sqrt{|\lambda|}$. Then for $T=1$ this is 
$O(e^{-\sqrt{|\lambda|}})$ and therefore the term
$$\int_T^\infty t^{s-d/2}e^{\lambda t}dt$$ 
is negligent.  This proves the asymptotic expansion in (\ref{ros}).

The right-hand side of  (\ref{ros}) is equal to 
$$
\sum_{n=0}^\infty \frac{a_n}{(-\lambda)^{s+n+1-d/2}}
\int_0^\infty u^{s+n-d/2}e^{-u}du =\sum_{n=0}^{\infty}
\frac{\Gamma(s+n+1-d/2)a_n}{(-\lambda)^{s+n+1-d/2}},
$$
and this completes the proof of the lemma.
\qed

\subsection{Agmon-Kannai expansion for derivatives of the resolvent}
In the notations of Theorem \ref{ak} let $H=\Delta$ be the Laplacian
on a $d$-dimensional Riemannian manifold $M$, and $\Delta_0$ be the 
operator obtained from the principal part of the Laplacian by freezing
its coefficients at a certain point $x\in M$. As before, 
$R_\lambda=(\Delta-\lambda)^{-1}$, $F_{\lambda}=(\Delta_0-\lambda)^{-1}$. 
\begin{theorem}
The following asymptotic expansion on the diagonal holds for the derivatives 
of the resolvent kernel of the Laplacian on a $d$-dimensional Riemannian
manifold $M$:
\begin{equation}
\label{der}
\frac{d^s}{d\lambda^s}R_\lambda(x,x)\sim 
\frac{1}{\sqrt{g}} \sum_{m=0}^{\infty}\frac{(m+s)!}{m!}
X_mF_{\lambda}^{m+s+1}, \quad s\ge d/2.
\end{equation}
\end{theorem}
\noindent {\bf Proof.}
Formally we have:
$$
\frac{d}{d\lambda}F_\lambda=\frac{d}{d\lambda}
\left(\frac{1}{\Delta_0-\lambda}\right)=\frac{1}{(\Delta_0-\lambda)^2}=
F_{\lambda}^2.
$$
This implies
$$
\frac{d^s}{d\lambda^s}F_\lambda^{m+1}=
\frac{(m+s)!}{m!}F_{\lambda}^{m+s+1}.
$$
Together with (\ref{s}) this completes the proof of the theorem. 
\qed

\medskip
Let us introduce the standard multi-index notations (see [H\"o]): 
if $\alpha=(\alpha_1,\dots,\alpha_d)$ is a multi-index, 
then $|\alpha|=\alpha_1+\cdots+\alpha_d$,
$\alpha!=\alpha_1!\cdots\alpha_d!$. 
For any vector $x=(x_1,\dots,x_d)$ we denote 
$x^{\alpha}=x_1^{\alpha_1}\cdots x_d^{\alpha_d}$ and
$$\frac{\partial^{\alpha}}{\partial x^{\alpha}}=
\frac{\partial^{\alpha_1}}{\partial x_1^{\alpha_1}}\cdots
\frac{\partial^{\alpha_d}}{\partial x_d^{\alpha_d}}.$$

We note that (\ref{der}) and (\ref{s}) are in fact asymptotic expansions
in the powers of $-\lambda$ as well as (\ref{res}). This is due to the
following formula (see [AK]):
\begin{multline}
\label{int}
\frac{\partial^{\gamma}}{\partial x^{\gamma}} 
F_\lambda^{m+s+1}(x,x)=\\
(-\lambda)^{\frac{d+|\gamma|}{2}-m-s-1}
\frac{(-1)^{\frac{|\gamma|}{2}}}{(2\pi)^d}\int\limits_{\Bbb R^d}
\frac{\xi^{\gamma}\,d\xi}
{(\Delta_0(\xi)^2+1)^{m+s+1}}, 
\end{multline}
where $\Delta_0(\xi)$ denotes the symbol of the operator
$\Delta_0$, $\gamma=(\gamma_1,\dots,\gamma_d)$ is a multi-index
and $\xi=(\xi_1,\dots,\xi_d)$.
\section{Heat invariants in normal coordinates}
\subsection{Computation of heat invariants}
Let $(x_1,\dots,x_d)$ be local coordinates on the 
Riemannian manifold $M$ 
such that the Riemannian metric at the origin $x=(0,\dots,0)\in M$ 
(in the sequel we simply write $x=0$), 
is Euclidean: $g_{ij}|_{x=0}=\delta_{ij}$ 
For convenience we may consider normal coordinates 
on $M$ centered at the point $x=0$ (see [GKM]).

\begin{theorem}
\label{main}
Let $M$ be a $d$-dimensional Riemannian manifold without boundary
and $(x_1,\dots,x_d)$ be normal coordinates on $M$ centered at the
point $x=0$. Then the heat invariants $a_n(x)$ at the point $x=0$ are 
equal to:
\begin{multline}
\label{q}
a_n(0)=
(4\pi)^{-\frac{d}{2}}(-1)^n 
\sum_{m=n}^{4n}\sum_{k=n}^m 
\frac{1}{k!\,2^{2m-2n}}\cdot \\
\cdot \sum_{|\alpha|=m-k}\,\,\sum_{|\beta|=k-n}
\frac{(2\alpha+2\beta)!}{\alpha!(\alpha+\beta)!(2\beta)!}
\left.\Delta^k(x^{2\beta})\right|_{x=0},
\end{multline}  
where 
$\alpha=(\alpha_1,\dots,\alpha_d)$, 
$\beta=(\beta_1,\dots,\beta_d)$ are multi-indices.
\end{theorem}
\noindent {\bf Proof of Theorem \ref{main}}
Since $(x_1,\dots,x_d)$ are normal coordinates centered at
$x=0$, the principal part of the Laplacian at this point
coincides with the Euclidean Laplacian, i.e.:
\begin{equation}
\label{d0}
\Delta_0=-\frac{\partial^2}{\partial x_1^2}-\cdots-
\frac{\partial^2}{\partial x_d^2}.
\end{equation}
Due to (\ref{res}), in order to compute the
coefficient $a_n(x)$ we have to collect all terms in the expansion
(\ref{der}) containing $(-\lambda)^{d/2-s-n-1}$.
From (\ref{int}) we have:
$$\frac{d+|\gamma|}{2}-m-s-1=\frac{d}{2}-s-n-1,$$
which implies $|\gamma|=2m-2n$ and in particular $m\ge n$.
As it was shown in [P1], estimates on the orders of operators $X_m$
(namely, Lemma 3.1 and Theorem 5.1 in [AK]) imply that $m\le 4n$. 

Note that due to (\ref{d0}) all indices $\gamma_1,\dots,\gamma_2$
should be even since otherwise the integral in (\ref{int}) will vanish.
Setting $\gamma=2\mu=(2\mu_1,\dots,2\mu_d)$ 
and taking into account that $|\mu|=m-n$ we compute this integral
(see [GR]):
$$
\int\limits_{\Bbb R^d}
\frac{\xi^{2\mu}\,d\xi}{(\xi^2+1)^{m+s+1}}=
\frac{\Gamma(\mu_1+\frac12)\Gamma(\mu_2+\frac12)\cdots\Gamma(\mu_d+\frac12)
\Gamma(s+n+1-\frac{d}{2})}{(m+s)!}
$$
Substituting this into (\ref{int}) and further on into (\ref{der})
we obtain due to (\ref{res}):
\begin{equation*}
a_n(x)=\left.\sum_{m=n}^{4n}\sum_{k=0}^m 
\frac{(-1)^{k+m-n}}{m!(2\pi)^d}\binom{m}{k}
\Delta^k\Delta_0^{m-k}\left(\sum_{|\mu|=m-n}\frac{x^{2\mu}}{(2\mu)!}
\prod_{i=1}^d\Gamma\left(\mu_i+\frac12\right)
\right)\right|_{x=0}
\end{equation*}
Note that $s$ has cancelled out as one could expect since heat 
invariants do not depend on $s$!

Now let us simplify this formula.
First notice that
\begin{equation}
\label{1}
\Delta_0^{m-k}=(-1)^{m-k}\sum_{|\beta|=m-k}\frac{(m-k)!}{\beta!}
\frac{\partial^{2\beta}}{\partial x^{2\beta}},
\end{equation}
where $\beta=(\beta_1,\cdots,\beta_d)$.
Using the well-known representation of the $\Gamma$-function
\begin{equation}
\label{gam}
\Gamma(k+1/2)=\frac{\sqrt{\pi}(2k)!}{4^k k!}
\end{equation}
we also obtain:
\begin{equation}
\label{2}
\prod_{i=1}^d\Gamma\left(\mu_i+\frac12\right)=
\frac{\pi^{d/2}}{2^{2m-2n}}\frac{(2\mu)!}{\mu!}.
\end{equation}
Let us substitute (\ref{1}) and (\ref{2}) into the above 
formula for the $a_n(x)$ and apply $\Delta_0^{m-k}$ to $x^{2\mu}$.
Introducing the new summation multi-index $\alpha=\mu-\beta$ 
and noticing that all terms for $k<n$ vanish we finally obtain
\begin{equation*}
a_n(0)=(4\pi)^{-\frac{d}{2}}(-1)^n \sum_{m=n}^{4n}\sum_{k=n}^m 
\frac{1}{k!\,2^{2m-2n}}\sum_{|\alpha|=m-k}\,\,\sum_{|\beta|=k-n}
\frac{(2\alpha+2\beta)!}{\alpha!(\alpha+\beta)!(2\beta)!}
\left.\Delta^k(x^{2\beta})\right|_{x=0}.
\end{equation*}  
This completes the proof of the theorem. \qed
\subsection{Remarks}
The proof of the Theorem \ref{main} is similar to the proofs of main
theorems in [P1] and [P2]. One may check that in the particular cases
of the $2$-dimensional Laplacian and the $1$-dimensional Schr\"odinger 
operator Theorem \ref{main} agrees with the results obtained 
in [P1] and [P2].
\section{Invariance and  combinatorial identities}
\subsection{Combinatorial identities}
Let us rewrite (\ref{q}) in the following way:
\begin{equation*}
(4\pi)^{-\frac{d}{2}}(-4)^n\sum_{k=n}^{4n}\frac{1}{k!}\left.
\Delta^k\left(\sum_{m=k}^{4n} 
\frac{1}{4^m}\sum_{|\beta|=k-n}\,\left(
\sum_{|\alpha|=m-k}
\frac{(2\alpha+2\beta)!\beta!}{\alpha!(\alpha+\beta)!(2\beta)!}
\right)
\frac{1}{\beta!}x^{2\beta}\right)\right|_{x=0}.
\end{equation*}  
Observe that due to the multinomial theorem
\begin{equation}
\label{mult}
\sum_{|\beta|=k-n}\frac{1}{\beta!}x^{2\beta}=\frac{1}{(k-n)!}
(x_1^2+\cdots +x_d^2)^{k-n}.
\end{equation}

Let us recall the following generalization of the binomial coefficients
(see [Er]). For real $z\in {\Bbb R}$ and $a\in {\Bbb N}$ set
\begin{multline}
\label{bc}
\binom{z}{a}=\binom{z}{z-a}=\\
\frac{\Gamma(z+1)}{\Gamma(a+1)\Gamma(z-a+1)}=
\frac{z(z-1)\cdots (z-a+1)}{a!}.
\end{multline}
We also set $\binom{z}{0}=\binom{z}{z}=1$.

Let us proceed  with the following simple combinatorial formula.
\begin{lemma} 
\label{hren}
\begin{equation}
\label{vspom2}
\sum_{a=0}^{u}\binom{z+a}{a}\binom{w+u-a}{u-a}=\binom{z+w+u+1}{z+w+1}.
\end{equation}
\end{lemma}
\noindent {\bf Proof.}
Using the method of generating functions (see [Rio]) we have:
$$
\sum_{a=0}\binom{z+a}{a}q^{2a}=\frac{1}{(1-q^2)^{z+1}}
$$
which implies
$$
\sum_{u=0}^{\infty}\sum_{a_1+a_2=u}\binom{z+a_1}{a_1}\binom{w+a_2}{a_2}
q^{2u}=\frac{1}{(1-q^2)^{z+w+2}}=
\sum_{u=0}^{\infty}\binom{z+w+u+1}{z+w+1}q^{2u}.
$$
This completes the proof of the lemma. \qed 

\medskip

Now we can prove our main combinatorial identity.
\begin{theorem}
\label{comb1}
Let $\alpha$, $\beta$ be multi-indices of dimension $d$ and let $|\beta|=v$.
Then
\begin{equation}
\label{expr}
\sum_{|\alpha|=u}
\frac{(2\alpha+2\beta)!\beta!}{\alpha!(\alpha+\beta)!(2\beta)!}=
4^u\,\binom{u+v-1+d/2}{u}.
\end{equation} 
\end{theorem}
\noindent{\bf Proof.}
We proceed by induction over $d$. For $d=1$ we have due to (\ref{gam}):
$$
\frac{(2u+2v)! v!}{u! (u+v)! (2v)!}=\frac{4^{u+v}\Gamma(u+v+1/2)
\sqrt{\pi}}{4^v \Gamma(v+1/2) u! \sqrt{\pi}}=4^u\binom{u+v-1/2}{u}. 
$$
and hence (\ref{expr}) is valid.

Suppose we have proved the formula (\ref{expr}) in all dimensions
less than some $d>1$. Let us prove it in the dimension $d$. 
Denote $\alpha_1=a$, $\beta_1=b$.
By induction we may rewrite the sum in (\ref{expr}) as
\begin{equation}
\label{vspom3}
\sum_{|\alpha|=u}
\frac{(2\alpha+2\beta)!\beta!}{\alpha!(\alpha+\beta)!(2\beta)!}=
\sum_{a=0}^u
\frac{(2a+2b)! b!}{a!(a+b)!(2b)!} \binom{u-a-1+l}{u-a}4^{u-a},
\end{equation}
where $l=v-b+\frac{d-1}{2}$.
On the other hand,
$$
\frac{1}{4^a}\frac{(2a+2b)!b!}{a!(a+b)!(2b)!}=\binom{a+b-1/2}{a}
$$
and hence (\ref{vspom3}) is equal to
$$
4^u \sum_{a=0}^u \binom{a+b-1/2}{a}\binom{u-a-1+l}{u-a}.$$ 
By Lemma (\ref{hren}) this equals to
$$4^u \binom{u+b+l-1/2}{u}=4^u \binom{u+v-1+d/2}{u},$$
which completes the proof of the theorem. \qed
\subsection{Proof of Theorem \ref{most}}
Set $u=m-k$, $v=k-n$.
Combining Theorem \ref{main}, Theorem \ref{comb1} and formula (\ref{mult})
we obtain the following reformulation of (\ref{q}):
\begin{multline}
\label{qq}
a_n(0)=\\
(4\pi)^{-d/2} (-4)^n \sum_{k=n}^{4n}
\left(\sum_{m=k}^{4n} \binom{m+\frac{d}{2}-n-1}{m-k}\right)
\frac{\Delta^k(|x|^{2k-2n})|_{x=0}}{k!\, 4^k},
\end{multline}
where $|x|^2=x_1^2+\cdots+x_d^2$.
Denote $i=m-k$, $j=k-n$.
By Lemma \ref{hren} the inner sum may be rewritten as
$$
\sum_{i=0}^{4n-k}\binom{i+d/2+k-n-1}{i}=\binom{3n+d/2}{j+d/2}.
$$
Therefore (\ref{qq}) is equal to
\begin{equation}
\label{qqq}
a_n(0)=(4\pi)^{-d/2} (-1)^n \sum_{j=0}^{3n}\binom{3n+d/2}{j+d/2}
\frac{\Delta^{j+n}(|x|^{2j})|_{x=0}}{4^j\, j! \,(j+n)!}.
\end{equation}

Consider the function $\rho_x(y)^2$ which is the square of the distance
between the points $x$ and $y$. 
In normal coordinates centered at the point $x=0$ it is given locally by
$$\rho_x(y)^2=\sum_{i,j=1}^d
g_{ij}(0)\,y_i y_j = y_1^2+\cdots+y_d^2=|y|^2$$ where 
$y=(y_1,..,y_d)$ (see [Du], p. 94).
Therefore we may rewrite formula (\ref{qqq}) in an invariant form, namely
\begin{equation*}
a_n(x)=(4\pi)^{-d/2} (-1)^n \sum_{j=0}^{3n}\binom{3n+d/2}{j+d/2}
\frac{\Delta_y^{j+n}(\rho_x(y)^{2j})|_{y=x}}{4^j\, j! \,(j+n)!},
\end{equation*}
where the subscript of the Laplacian means that the operator is acting on
functions in $y$-variable.
This completes the proof of Theorem \ref{most}. \qed

\section{Application to computation of the Korteweg-de Vries hierarchy}
\subsection{Asymptotics of Schr\"odinger operator and KdV hierarchy}
In [P2] we have applied the Agmon-Kannai method to computation of the 
Korteweg-de Vries hierarchy (see [NMPZ]).
Let us briefly recall the setting of the problem.

Consider the $1$-dimensional Schr\"odinger operator:
\begin{equation*}
L=\frac{\partial^2}{\partial x^2}+U(x).
\end{equation*}
Its heat kernel $H(t,x,y)$ is  
the fundamental solution of the heat equation 
$$\left(\frac{\partial}{\partial t} - L\right)f=0.$$
It has the following  asymptotic representation on the diagonal 
as $t\to 0+$: 
\begin{equation*}
H(t,x,x)\sim\frac{1}{\sqrt{4\pi t}}\sum_{n=0}^{\infty}h_n[U]t^n,
\end{equation*}
where $h_n[U]$ are some polynomials in $U(x)$ and 
its derivatives. 

The KdV hierarchy is defined by (see [AvSc]):
\begin{equation}
\label{kdvh}
\frac{\partial U}{\partial t}=\frac{\partial}{\partial x}G_n[U],
\end{equation}
where
\begin{equation*}
G_n[U]=
\frac{(2n)!}{2\cdot n!}h_n[U], 
\quad n\in {\Bbb N}.
\end{equation*}
Set $U_0=U$, $U_n=\partial^n U/\partial x^n$, $n\in {\Bbb N}$,
where $U_n$, $n\ge 0$ are formal variables. 
The sequence of polynomials $G_n[U]=G_n[U_0,U_1,U_2,\dots]$ 
starts with (see [AvSc]): 
\begin{equation*}
G_1[U]=U_0, \,\, G_2[U]=U_2+3U_0^2, \,\, 
G_3[u]=U_4+10U_0U_2+5U_1^2+10U_0^3, \, \dots
\end{equation*}
In particular, substituting $G_2[U]$ into (\ref{kdvh}) 
we obtain the familiar Korteweg-de Vries  equation (see [NMPZ]):
\begin{equation*}
\frac{\partial U}{\partial t}
=\frac{\partial^3 U}{\partial x^3}
+6U\frac{\partial U}{\partial x}
\end{equation*}
\subsection{Computation of the KdV hierarchy}
In [P2] we have presented explicit formulas for the KdV hierarchy 
(we refer to [P2] for the history of this question). Theorem \ref{most}
allows to simplify the results of [P2].

\begin{theorem}
\label{check}
The KdV hierarchy is  given by:
$$
G_n[U]=\frac{(2n)!}{2\cdot n!}\sum_{j=0}^{n}
\binom{n+\frac{1}{2}}{j+\frac{1}{2}}
\frac{(-1)^j}{4^j \, j! \, (j+n)!}P_{nj}[U],
\eqno{(\ref{check})}
$$
where the polynomial $P_{nj}[U]$ is obtained from 
$\left.L^{j+n}(x^{2j})\right|_{x=0}$ by a formal change of variables:
$U_i(0)\to U_i$, $i=0,..,2n+2j-2$.
\end{theorem}

This expression can be completely expanded due to a formula 
for the powers of the Schr\"odinger operator ([Rid]).

\begin{theorem}
\label{KdV1}
The  polynomials $G_n[U]$, $n\in {\Bbb N}$ are equal to:
$$
G_n[U]=
\frac{(2n)!}{2\cdot n!}\sum_{j=0}^{n}
\binom{n+\frac{1}{2}}{j+\frac{1}{2}}
\frac{(-1)^j (2j)!}{4^j \, j! \, (j+n)!}
\sum_{p=1}^{j+n}\mskip-2\thinmuskip
\sum_{k_1,\dots,k_p\atop k_1+\cdots +k_p=2(n-p)}\mskip-13\thinmuskip
C_{k_1,\dots,k_p}U_{k_1}\cdots U_{k_p},
$$
where 
\begin{equation*}
C_{k_1,\dots,k_p}=\sum\begin{Sb}
0\le\l_0\le l_1\le \cdots\le l_{p-1}=j+n-p\\ 
2l_i\ge k_1+\cdots+k_{i+1},\,\,i=0,\dots,p-1.\end{Sb}
\binom{2l_0}{k_1}\binom{2l_1-k_1}{k_2}\cdots
\binom{2l_{p-1}-k_1-\cdots-k_{p-1}}{k_p}.
\end{equation*}
\end{theorem}      

\subsection*{Remark}
Formula (\ref{check}) was checked using Mathematica ([Wo]) 
and for $1\le n\le5$ the results agreed with the 
already known ones (cf. [GD]).
\subsection*{Acknowledgments}
This paper is a part of my Ph.D. research at the Department of Mathematics of 
the Weizmann Institute of Science. I am very grateful to my Ph.D. advisor 
Yakar Kannai for his constant  help and support. 

The author is indebted
to Leonid Polterovich, Amitai Regev and Mikhail Solomyak for helpful 
discussions. I would like to thank Ivan Avramidi, Isaac Chavel,  
Peter Gilkey and Steven Rosenberg for their remarks on the preliminary 
versions of this paper. 
I am also grateful to Sergei Novikov for a useful discussion on the KdV
hierarchy and its computation.

Theorem \ref{comb1} was first established with the help of a Mathematica
implementation of the Wilf-Zeilberger algorithm (see [PWZ]). 
I am thankful to Marko Petkov\v{s}ek and Doron Zeilberger for their help with 
this matter.

\section*{References}

\bigskip

\noindent [AK] S. Agmon, Y. Kannai, On the asymptotic behavior of spectral 
functions and resolvent kernels of elliptic operators, Israel J. Math. 5
(1967), 1-30.

\smallskip

\noindent [Av] I.G. Avramidi, A covariant technique for the calculation of the 
one-loop effective action, Nucl. Phys. B 355 (1991) 712-754; 
Erratum: Nucl. Phys. B 509 (1998) 557-558.

\smallskip

\noindent [AvB] I.G. Avramidi, T. Branson, Heat kernel asymptotics of operators
with non-Laplace principal part, math-ph/9905001, (1999), 1-43.

\smallskip

\noindent [AvSc] I.G. Avramidi and R. Schimming, A new explicit expression
for the Korteweg - de Vries hierarchy, solv-int/9710009,  (1997), 1-17.

\smallskip
 
\noindent [Be] M. Berger, Geometry of the spectrum, Proc. Symp. Pure Math. 27 
(1975), 129-152.

\smallskip

\noindent [BGM] M. Berger, P. Gauduchon, E. Mazet, Le spectre d'une vari\`et\`e
Riemannienne, Lecture Notes in Math. 194, Springer-Verlag, 1971.

\smallskip

\noindent [BG\O] T. Branson, P. Gilkey, B. \O rsted, Leading terms in the heat 
invariants, Proc. Am. Math. Soc. 109 (1990), 437-450. 

\smallskip

\noindent [Ch] I. Chavel, Eigenvalues in Riemannian geometry, Academic Press, 
1984.

\smallskip

\noindent [Da] E. B. Davies, Heat kernels and spectral theory.
Cambridge University Press, 1989.
 
\smallskip

\noindent [Du] G.F.D. Duff, Partial differential equations, University
of Toronto Press, 1956.

\smallskip

\noindent [Er] A. Erd\'elyi et. al., Higher transcendental functions, vol. 1,
McGraw-Hill, 1953.

\smallskip

\noindent [F] S.A. Fulling ed., Heat Kernel Techniques and Quantum Gravity,
Discourses in Math. and its Appl., No. 4, Texas A\&M Univ., 1995.

\smallskip

\noindent [GD] I.M. Gelfand, L.A. Dikii, Asymptotic behaviour of the
resolvent of Sturm-Liouville equation and the algebra of the 
Korteweg-de Vries equations, Russian Math. Surveys, 30:5 (1975),
77-113.

\smallskip

\noindent [G1] P. Gilkey, The spectral geometry of a Riemannian manifold, J. 
Diff. Geom. 10 (1975), 601-618.   

\smallskip

\noindent [G2]  P. Gilkey, The index theorem and the heat equation, Math. 
Lect.Series,  Publish or Perish, 1974.

\smallskip

\noindent [G3] P. Gilkey, Heat equation asymptotics, Proc. Symp. Pure Math. 
54 (1993), 317-326.

\smallskip

\noindent [GR] I.S. Gradshtein, I.M. Ryzhik, Table of integrals, series and
products, Academic Press,  1980.

\smallskip

\noindent [GKM] D. Gromoll, W. Klingenberg, 
W. Meyer, Riemannsche Geometrie im Gros\-sen, Springer-Verlag, 1968.

\smallskip

\noindent [H\"o] L. H\"ormander, 
Linear partial differential operators, Springer-Verlag, 1969.

\smallskip

\noindent [MS]  H.P. McKean, Jr., I.M. Singer, 
Curvature and the eigenvalues of the 
Laplacian, J. Diff. Geom. 1 (1967), 43-69.

\smallskip

\noindent [MP] S. Minakshisundaram, A. Pleijel, 
Some properties of the eigenfunctions
of the Laplace-operator on Riemannian manifolds, Canadian J. Math. 1 (1949),
242-256.

\smallskip

\noindent [NMPZ] S.P. Novikov, S.V. Manakov, L.P. Pitaevskii, V.E. Zakharov,
Theory of solitons: the inverse scattering method, 
Consultants Bureau, 1984. 

\smallskip

\noindent [OPS] B. Osgood, R. Phillips, P. Sarnak, Compact isospectral sets of
surfaces, J. Funct. Anal. 80 (1988), 212-234.

\smallskip

\noindent [PWZ] M. Petkov\v{s}ek, 
H. Wilf, D. Zeilberger, A=B, A K Peters, 1996.

\smallskip

\noindent [P1] I. Polterovich, A commutator method for computation of heat
invariants,  to appear in Indagationes Mathematicae, 1-11.

\smallskip

\noindent [P2] I. Polterovich, From Agmon-Kannai expansion to Korteweg-de Vries
hierarchy, to appear  in Letters in Mathematical Physics, 1-7.

\smallskip

\noindent [Rid] S.Z. Rida, Explicit formulae for the powers of a 
Schr\"odinger-like
ordinary differential operator, J. Phys. A: Math. Gen. 31 (1998), 5577-5583.

\smallskip

\noindent [Rio] J. Riordan, An introduction to combinatorial analysis,
John Wiley \& Sons, Inc., 1958. 

\smallskip

\noindent [Ro] S. Rosenberg, The Laplacian on a Riemannian manifold, Cambridge
University Press, 1997.

\smallskip

\noindent 
[Se] R. Seeley, Complex powers of an elliptic operator, Proc. Symp. Pure
Math. 10 (1967), 288-307.

\smallskip

\noindent 
[Sa] T. Sakai, On the eigenvalues of the Laplacian and curvature of Riemannian
manifold, T\^ohuku Math. J. 23 (1971), 585-603.

\smallskip

\smallskip

\noindent [vdV] A.E.M. van de Ven, Index--free heat kernel coefficients,
hep-th/9708152, (1997), 1--38.

\smallskip

\noindent [Wo] S.Wolfram, Mathematica: a system for doing mathematics 
by computer, Addison -- Wesley, 1991.

\smallskip

\noindent [Xu] C. Xu, Heat kernels and geometric invariants I, Bull. Sc. Math.
     117 (1993), 287-312.

\end{document}